\documentclass[12pt]{amsart}
\usepackage[all]{xy}
\usepackage{amssymb}
\usepackage{amsthm}
\usepackage{amsmath}
\usepackage{amscd,enumitem}
\usepackage{verbatim}
\usepackage{eurosym}
\usepackage{graphicx}
\usepackage{dsfont}
\usepackage{stmaryrd}

\usepackage{float}
\usepackage{color}
\usepackage{longtable}
\usepackage{dcolumn}
\usepackage[mathscr]{eucal}
\usepackage[all]{xy}
\usepackage{hyperref}
\usepackage[title]{appendix}
\usepackage[usenames,dvipsnames]{xcolor}
\usepackage{bbm}
\usepackage[textheight=8.5in, textwidth=6.5in]{geometry}
\usepackage{multirow}
\usepackage{caption}
\newtheorem*{thm*}{Theorem}
\newtheorem*{conj*}{Conjecture}
\newtheorem{theorem}{Theorem}

\newtheorem*{remark}{Remark}

\newtheorem{thm}{Theorem}[section]

\newtheorem{prop}[thm]{Proposition}

\newtheorem*{example}{Example}

\newcommand{\Z}{\mathbb{Z}}
\newcommand{\Q}{\mathbb{Q}}

\newcommand{\HH}{\mathbb{H}}
\newcommand{\DD}{\mathbb{D}}
\newcommand{\PP}{\mathbb{P}}

\newcommand{\SL}{\operatorname{SL}}

\numberwithin{equation}{section}

\begin{document}
\title[Elliptic expansion Moonshine Conjecture]
{Proof of the elliptic expansion Moonshine Conjecture of
 C\u{a}ld\u{a}raru, He, and Huang} \author{Letong Hong, Michael H. Mertens, Ken Ono, \and Shengtong Zhang}

 \address{Department of Mathematics, Massachusetts Institute of Technology,
Cambridge, MA 02139}
\email{clhong@mit.edu}

\address{Department of Mathematical Sciences, University of Liverpool, Liverpool L69 7ZL}
\email{m.h.mertens@liverpool.ac.uk}

\address{Department of Mathematics, University of Virginia, Charlottesville, VA 22904}
 \email{ken.ono691@virginia.edu}
 
 \address{Department of Mathematics, Massachusetts Institute of Technology,
Cambridge, MA 02139}
\email{stzh1555@mit.edu}


\begin{abstract} Using predictions in mirror symmetry, C\u{a}ld\u{a}raru, He, and Huang recently formulated a ``Moonshine Conjecture at Landau-Ginzburg points''  \cite{CHH} for Klein's modular $j$-function at $j=0$ and $j=1728.$ The conjecture asserts that the $j$-function, when specialized at
specific flat coordinates on the moduli spaces of versal deformations of the corresponding CM elliptic curves,
yields simple rational functions. We prove this conjecture, and show that these rational functions arise from classical $ _2F_1$-hypergeometric inversion formulae for the $j$-function.
\end{abstract}
\maketitle
\section{Introduction and statement of results}

The {\it Monstrous Moonshine Conjecture} \cite{CN}, famously proved by Borcherds \cite{Borcherds},
  offers a  surprising relationship between the  Monster $\mathbb{M},$ the  largest sporadic simple group, and generators for the modular function fields of specific genus 0 modular curves.  The conjecture sprouted from the observation that the first few coefficients of Klein's function
\begin{equation}\label{J-function}
j(\tau)-744=q^{-1}+196884q+21493760q^2+O\left(q^3\right),
\end{equation}
the Hauptmodul for $\mathrm{SL}_2(\mathbb{Z})$ with $q:=e^{2\pi i\tau}$ and $\tau$ in the upper-half of the complex plane, are simple sums of the dimensions of irreducible representations of $\mathbb{M}$.  
Monstrous Moonshine  \cite{Borcherds, CN, Thompson1, Thompson2} offers
 a graded, infinite-dimensional $\mathbb{M}$-module $$V^\natural=\bigoplus_{n\geq-1} V^\natural(n)$$ whose graded dimensions are the coefficients of $j(\tau)-744$.   Moreover, 
 for each $g\in\mathbb{M}$, it provides a genus zero subgroup $\Gamma_g\subseteq\mathrm{SL}_2(\mathbb{R})$ for which 
$$
T_g(\tau):=\sum_{n\geq-1}\mathrm{Tr}\left(g\vert V^\natural(n)\right)q^n,
$$
is the Hauptmodul for $\Gamma_g$.

Many further examples of ``moonshine'' have been obtained (for example, see \cite{CDH,Gannon,DGO15,DMO17}), and these advances have often been inspired by physical considerations (e.g. black hole entropies, quantum gravity, etc.). These works illustrate that there are many more types of moonshine than mathematicians and physicists initially thought.  For a survey of moonshine and its applications to mathematics and physics, the reader may consult \cite{DGO1} or \cite{GannonBook}.

In a recent paper \cite{CHH}, C\u{a}ld\u{a}raru, He, and Huang discovered a new remarkable phenomenon which they refer to as ``Moonshine at Landau-Ginzburg points\footnote{Here  the term ``moonshine'' refers to the {\it foolish idea} that this combinatorial geometric specialization of the $j$-function reduces to a simple rational function.}.''  Instead of making use of the Fourier expansion of $j(\tau)$ at the cusp infinity, their conjecture involves the elliptic expansions at $\tau_{*}=\rho:=e^{\pi i/3}$ and $\tau_*=i.$   Combining predictions from mirror symmetry with  recent computations of higher genus Fan-Jarvis-Ruan-Witten (FJRW) invariants by Li, Shen, and Zhou \cite{LSZ}, they conjectured that the elliptic expansion at $\rho$ (resp. $i$),  specialized at the normalized flat coordinate of the corresponding moduli space of versal deformations of elliptic curves, is the classical rational function defining the hypergeometric inversion formulae for the $j$-function (see Proposition~\ref{PiAGM_rho} and Proposition~\ref{PiAGM_i}).

To recall their conjectures, let $\HH$ denote the upper-half of the complex plane, and let $\DD$ be the unit disk.  To define the Taylor expansions of the $j$-function at $\tau_{*}\in \{ \rho, i\},$ we employ the
uniformizing map  $S_{\tau_*}: \HH \to \DD$  and its inverse
$S_{\tau*}^{-1}:\DD\to \HH$ defined by
$$
S_{\tau_*}(\tau):=\frac{\tau-\tau_{*}}{\tau-\overline{\tau}_{*}} \ \ \ \ \  {\text {\rm and}}\ \ \ \ \ 
S_{\tau_*}^{-1}(w):=\frac{\tau_{*}-\overline{\tau}_{*}w}{1-w}.
$$
Finally, we define the renormalized inverse $s_{\tau_*}^{-1}$  by
\begin{equation}
s_{\tau_*}^{-1}(w):= S_{\tau_*}^{-1}\left( \frac{w}{2\pi \Omega_{\tau_*}^2}\right) ,
\end{equation}
where $\Omega_{\tau_*}$ are the standard Chowla-Selberg periods (for example, see
(96) of Section~6.3 of \cite{Zagier123})
\begin{equation}
\Omega_{\tau_*}:=\begin{cases} \frac{1}{\sqrt{6\pi}}\cdot \left(\frac{\Gamma(1/3)}{\Gamma(2/3)}\right)^{\frac{3}{2}}\ \ \ \ 
\ &{\text {\rm if $\tau_*=\rho$}},\\ 
\frac{1}{\sqrt{8\pi}}\cdot \left(
\frac{\Gamma(1/4)}{\Gamma(3/4)}\right) \ \ \ \ &{\text {\rm if $\tau_*=i$.}}
\end{cases}
\end{equation}
As power series in $w$, we have the following expansions
\begin{equation}\label{eqexpansion}
\begin{split}
j\left(s_{\rho}^{-1}(w)\right)=13824w^3-39744w^6+\frac{1920024}{35}w^9-\dots=:\sum_{n=0}^{\infty} b_\rho(n)w^n,\\
j\left(s_{i}^{-1}(w)\right)=1728+20736w^2+105984w^4+\frac{1594112}{5}w^6+\dots=:\sum_{n=0}^{\infty} b_i(n)w^n.
\end{split}
\end{equation}
It is not difficult to find a recursive formula for the coefficients in \eqref{eqexpansion} by following the method of Proposition 28 of \cite{Zagier123}. For example, we find that $b_\rho(n)=-1728 \cdot 6^n p_n(0)/n!$, where the $p_n(t)\in\mathbb{Q}[t]$ form a sequence of polynomials $p_0(t)=t^3, p_1(t)=-t^2,\dots$, where for $n\geq 2$ we have
$$
 p_{n}(t)=\frac{t^3-1}{3}p_{n-1}'(t)-\frac{n+5}{6}t^2p_{n-1}(t)-\frac{(n-1)(n-2)}{144}tp_{n-2}(t).
$$
In particular, it follows that the coefficients $b_{\rho}(n)$ and $b_i(n)$ are all rational.

The conjecture of C\u{a}ld\u{a}raru, He, and Huang gives a striking formula
for the $t$-series
\begin{equation}\label{conjectures}
j\left(s_{\tau_{*}}^{-1}(c_{\tau_*}(t)\right)=\begin{cases} 13824t^3-46656t^6+99144t^9-171315t^{12}+\dots  \ \ \ \ \ \
&{\text {\rm if $\tau_{*}=\rho,$}}  \\
1728+20736t^2+147456t^4+851968t^6+\dots
\ \ \ \ \ \ &{\text {\rm if $ \tau_*=i,$}}
\end{cases}
\end{equation}
where $w=c_{\tau_{*}}(t)\in \Q\llbracket t\rrbracket,$
is the normalized flat coordinate (see Proposition~\ref{Flat})  of a corresponding moduli space of versal deformations 
of specific CM elliptic curves.
In the case of $\rho,$ the Fermat cubic (i.e. $j=0$) in $\PP^2$  is defined by
$$x^3+y^3+z^3=0,
$$
and the flat coordinate $c_{\rho}(t)$ corresponds to the moduli space of its versal deformations, the Hesse pencil of elliptic curves
\begin{equation}\label{HessePencilEC}
E_t:\ \  x^3+y^3+z^3+3txyz=0.
\end{equation}
For $j=1728$, the flat coordinate $c_{i}(t)$ arises similarly from the quartic in $\PP_{1,1,2}^2$ defined by
\begin{equation}\label{quartic}
x^4+y^4+z^2=0.
\end{equation}
Assuming this notation, we give a reformulation of their conjecture.

\begin{conj*}[C\u{a}ld\u{a}raru, He, and Huang \cite{CHH}, 2021]
The following are true.

\begin{enumerate}
\item If $|t|<1,$ then 

$$
j\left(\frac{s_{\rho}^{-1}(c_{\rho}(t)) + 1}{3}\right) = 27t^3\left(\frac{8-t^3}{1+t^3}\right)^3.
$$

\item  If $|t|<1/2,$ then
$$
j\left(s_{i}^{-1}(c_{i}(t))\right)= 
64\cdot  \frac{(3+4t^2)^3}{(1-4t^2)^2}.
$$\end{enumerate}

\end{conj*}

\medskip
\noindent
{\bf Three remarks.}

\smallskip
\noindent
(1)  The conjecture was originally formulated as identities between formal power series. Here we reformulate\footnote{We also correct a typographical error in part (2) of their conjecture.}  the conjecture in terms of analytic functions near $t=0.$  In particular, this reformulation introduces the denominator $3$ in the argument of $j$ in part (1) of the conjecture. 

\smallskip

\smallskip
\noindent
(2) The rational functions in the conjecture appear naturally as $j$-invariants.
For the Hesse pencil of elliptic curves $E_t$ in (\ref{HessePencilEC}), we have that $j(E_t)=27t^3\left(\frac{8-t^3}{1+t^3}\right)^3$. Similarly, for 
$$\mathcal{E}_{\lambda}: \ \ y^2=(x-1)\left (x^2-\frac{1}{\lambda+1}\right),
$$
we have $j(\mathcal{E}_{4t^2-1})=64\cdot  \frac{(3+4t^2)^3}{(1-4t^2)^2}.$  The symmetric square of $\mathcal{E}_{\lambda}$ is prominent in that it defines a parameterized family of $K3$ surfaces with generic Picard number 19 studied by Ahlgren, Penniston, and the third author \cite{AOP}.
\smallskip

\noindent
(3) The reader might also recognize these rational functions as they arise in classical inversion formulae
for the $j$-function (see Proposition~\ref{PiAGM_rho} and ~\ref{PiAGM_i}) involving the ${_2F_1(1/3, 2/3, 1;z)}$ and ${_2F_1(1/2, 1/2; 1;z)}$ hypergeometric functions.
 There are further inversion formulae, referred to as Ramanujan's theories of elliptic functions to alternative bases, involving the hypergeometric functions
$_2F_1(1/4, 3/4;1;z)$ and $ _2F_1(1/6, 5/6; 1; z)$
(for example, see \cite{BBG}).  Hence, it is natural to expect similar results  for corresponding families of CM elliptic curves.

\smallskip

In this note we prove the following theorem.

\begin{theorem}\label{main}
The C\u{a}ld\u{a}raru, He, and Huang Conjecture is true.
\end{theorem}

\begin{remark}  When this paper was submitted for publication, the authors were
 unaware of the 2018 paper \cite{SZ} by Shen and Zhou (see Section 3.3.3 of \cite{SZ}), which was written even before 
  C\u{a}ld\u{a}raru, He, and Huang
  formulated their two conjectures. Shen and Zhou explicitly note and derive the expressions
  for the expansions at $\rho$ in terms of the hypergeometric functions in this paper. Although they do not explicitly state or claim Conjecture (1) as a result, it is fair to say that their earlier work proves Conjecture (1).
  Therefore, the main contribution in this work is the proof of Conjecture (2).
\end{remark}

To prove Theorem~\ref{main}, we require a few important facts.
The most prominent ingredients are the descriptions of the normalized flat coordinates $c_{\tau_{*}}(t)$ obtained by Tu \cite{Tu19} and  C\u{a}ld\u{a}raru, He, and Huang
\cite{CHH}.  It turns out that the $c_{\tau_{*}}(t)$  at $\tau_* = \rho$ and $\tau_* = i$ possess simple descriptions as
 $ _2F_1$ hypergeometric functions. The proof of Theorem~\ref{main} makes use of 
 suitable transformation laws for such hypergeometric functions, some combinatorial considerations, and classical inversion formulae for the $j$-function in terms of hypergeometric functions. The ingredients for the proof are given in Section 2, and the proof is given in Section 3. Finally, in Section 4 we offer some examples of Theorem 1.

\section*{Acknowledgements} \noindent
The authors thank Andrei C\u{a}ld\u{a}raru for useful discussions.  The third
  author thanks  the Thomas Jefferson Fund and the NSF
(DMS-2002265 and DMS-2055118) for their generous support, and the support he received from  the Kavli Institute grant NSF PHY-1748958. Finally, the authors thank the referee for bringing their attention to the work of Shen and Zhou \cite{SZ} that proved Conjecture (1) in 2018.

\section{Nuts and bolts}\label{NutsAndBolts}

Here we recall the necessary ingredients for the proof of Theorem~\ref{main}.

\subsection{Explicit formulas for $c_{\tau_*}(t)$}

The following proposition gives the flat coordinates $c_{\tau_{*}}(t)$ obtained by Tu (see Section~4 of \cite{Tu19}) and  C\u{a}ld\u{a}raru, He, and Huang
(see Section~1.5 of \cite{CHH}).  For completeness,  we recall the multi-factorial notation.  For positive integers $M,n,r$ with $r<M$ we define
$$
(Mn-r)!!...!:=\prod_{m=1}^n(Mm-r),
$$
where we write $M$ exclamation marks on the left-hand side.

\begin{prop}\label{Flat} The following identities are true.
\begin{enumerate}
\item We have that $c_{\rho}(t)=h_{\rho}(t)/g_{\rho}(t),$ where
\begin{displaymath}
\begin{split}
h_{\rho}(t)&:=\sum_{n=0}^{\infty} (-1)^n \frac{\left((3n-1)!!!\right)^3}{(3n+1)!} t^{3n+1},\\
g_{\rho}(t)&:=\sum_{n=0}^{\infty} (-1)^n\frac{\left((3n-2)!!!\right)^3}{(3n)!} t^{3n}.
\end{split}
\end{displaymath}
\item We have that $c_i(t)=h_i(t)/g_i(t),$ where
\begin{displaymath}
\begin{split}
h_i(t):=\sum_{n=0}^{\infty}\frac{\left((4n-1)!!!!\right)^2}{(2n+1)!}t^{2n+1},\\
g_i(t):=\sum_{n=0}^{\infty}\frac{\left((4n-3)!!!!\right)^2}{(2n)!}t^{2n}.
\end{split}
\end{displaymath}
\end{enumerate}
\end{prop}

\subsection{Hypergeometric functions}

To prove Theorem~\ref{main}, we require some classical facts about $_2F_1$ hypergeometric functions, which are defined by
\begin{equation}
_2F_1(a, b; c; z):=\sum_{n=0}^{\infty}\frac{(a)_n (b)_n}{(c)_n}\cdot \frac{z^n}{n!},
\end{equation}
where the usual Pochhammer symbol is defined by
$$
(q)_n:=\begin{cases} 1 \ \ \ \ &{\text {\rm if $n=0,$}}\\
q(q+1)\cdots (q+n-1) \ \ \ \ &{\text {\rm if $n>0.$}}
\end{cases}
$$
We require the following two classical hypergeometric transformation laws.

\begin{prop}[Equation {15.10.33} of \cite{NIST}]\label{hyper1}
Assuming that the quantities are well-defined, for $|\arg(z)|<\pi$ we have: 
\begin{displaymath}
\begin{split}
_2F_1&(a,b;c;z)=\frac{\Gamma(1-b)\Gamma(c)}{\Gamma(a-b+1)\Gamma(c-a)}\left(\frac{1}{z}\right)^a \cdot  \ _2F_1\left (a-c+1,a;a-b+1;\frac{1}{z}\right)\\ &\ \ \ +
\frac{\Gamma(1-b)\Gamma(c)}{\Gamma(a)\Gamma(c-a-b+1)}\left(1-\frac{1}{z}\right)^{c-a-b}
\left(-\frac{1}{z}\right)^b \cdot  \ _2F_1\left(c-a,1-a; c-a-b+1;1-\frac{1}{z}\right).
\end{split}
\end{displaymath}
\end{prop}
\medskip
\noindent
{\bf Remark.}
The identity in Proposition \ref{hyper1} is obtained from Equation 15.10.33 of \cite{NIST}
by combining 15.10.14--15.

\begin{prop}[Equations 15.8.27--28 of \cite{NIST}]\label{hyper2}
Assuming that the quantities are well-defined, for $|\arg(z)|<\pi$ we have:\vspace{14pt}\\
(1)\vspace{-28pt}
\begin{displaymath}\hspace{-46pt}
\begin{split}
&\frac{2\Gamma(1/2)\Gamma(a+b+\frac12)}{\Gamma(a+\frac12)\Gamma(b+\frac12)} \ _2F_1(a,b;1/2;z)\\&\qquad ={_2F_1\left(2a,2b;a+b+\frac12;\frac{1-\sqrt{z}}{2}\right)}+ {_2F_1\left(2a,2b;a+b+\frac12;\frac{1+\sqrt{z}}{2}\right)}.
\end{split}
\end{displaymath}\\
(2)\vspace{-28pt}
\begin{displaymath}
\begin{split}
\qquad&\frac{2\sqrt{z}\Gamma(-1/2)\Gamma(a+b-\frac12)}{\Gamma(a-\frac12)\Gamma(b-\frac12)} \ _2F_1(a,b;3/2;z)\\&\quad = {_2F_1\left(2a-1,2b-1;a+b-\frac12;\frac{1-\sqrt{z}}{2}\right)}- {_2F_1\left(2a-1,2b-1;a+b-\frac12;\frac{1+\sqrt{z}}{2}\right)}.
\end{split}
\end{displaymath}
\end{prop}

\subsection{Hypergeometric functions and $j$-invariants}
Here we recall two classical inversion formulae for $j$-invariants of complex elliptic curves phrased in terms of $ _2F_1$ hypergeometric functions.

\begin{prop}[Theorem 4.4--4.5 of \cite{BBG}, or Eq. (2.8) of \cite{BC}]\label{PiAGM_rho}
 If $\tau\in \HH$ and $\gamma$ satisfies $$
\tau=\frac{i}{\sqrt{3}}\cdot \frac{ _2F_1\left(1/3, 2/3;1; 1-\gamma\right)}{ _2F_1 \left(1/3, 2/3; 1; \gamma\right)}
,$$ then we have $$j(\tau)=\frac{27(1+8\gamma)^3}{\gamma(1-\gamma)^3}.$$

\end{prop}

\begin{prop}[Corollary 5.17 of \cite{BBG}, or p. 431 of \cite{BC}]\label{PiAGM_i}
If $\tau\in \HH$ and $\lambda$ satisfies $$\tau=i\cdot \frac{_2F_1(1/2, 1/2; 1; 1-\lambda)}{_2F_1(1/2, 1/2; 1; \lambda)},$$ then we have $$j(\tau)=\frac{256(1-\lambda+\lambda^2)^3}{\lambda^2(1-\lambda)^2}.$$
\end{prop}

\section{Proof of the conjecture}

We now combine the facts from the previous section to prove Theorem~\ref{main}.  

\smallskip
\noindent
{\bf Case of $\rho$:}  
The case where $t=0$ follows by direct calculation. Therefore, we may assume
  that $0<|t|<1,$ which ensures
 the convergence of the functions in this proof.
 Now we note that
\begin{displaymath}
\begin{split}
_2F_1(1/3, 1/3; 2/3; -t^3)&=\sum_{n=0}^{\infty}\frac{(-t^3)^n}{n! (2/3)(5/3)\cdots
((3n-1)/3)}\cdot \left(\frac{1}{3}\cdot \frac{4}{3}\cdots \frac{3n-2}{3}\right)^2\\
&=\sum_{n=0}^{\infty}(-1)^n\frac{t^{3n}}{(3n)!}\prod_{m=1}^n (3m-2)^3\\
&=g_{\rho}(t).
\end{split}
\end{displaymath}
For convenience, we first assume that $0\leq |\arg(t)|< \pi/3.$ Applying Proposition~\ref{hyper1} with $a=1/3, b=1/3, c=2/3,$ and $z=-t^3,$ we find (using the standard branch of the cube root) that
\begin{equation}\label{rho1}
g_{\rho}(t)
=\frac{\Gamma(2/3)^2}{\Gamma(1/3)} t^{-1}\left( {_2F_1\left(1/3, 2/3; 1+t^{-3}\right)}
+\rho\cdot {_2F_1\left(1/3, 2/3; 1; -t^{-3}\right)}\right).
\end{equation}
Similarly, we have that
\begin{displaymath}
\begin{split}
t \cdot {_2F_1\left( 2/3, 2/3; 4/3; -t^3\right)}&=\sum_{n=0}^{\infty}
t\cdot \frac{(-t^3)^n}{n! (4/3)\cdots ((3n+1)/3)}\cdot \left(
\frac{2}{3}\cdots \frac{3n-1}{3}\right)^2\\
&=\sum_{n=0}^{\infty} (-1)^n \frac{t^{3n+1}}{(3n+1)!}\cdot \prod_{m=1}^n (3m-1)^3\\
&=h_{\rho}(t).
\end{split}
\end{displaymath}
Applying Proposition~\ref{hyper1} with $a=2/3, b=2/3, c=4/3,$ and $z=-t^3$ gives
\begin{equation}\label{rho2}
h_{\rho}(t)=\frac{\Gamma(1 / 3)^2}{3\Gamma(2 / 3)}t^{-1} \left({_2F_1(1/3, 2/3;1; 1 + t^{-3}}) -\overline{\rho}\cdot {_2F_1}(1/3, 2/3;1; -t^{-3})\right).
\end{equation}
Therefore, by combining (\ref{rho1}) and (\ref{rho2}) with Proposition~\ref{Flat}(1), we find 
 that
\begin{equation}\label{crhoformula}
c_{\rho}(t) = \frac{\Gamma(1 / 3)^3}{3\Gamma(2 / 3)^3}\  \cdot \frac{{_2F_1}(1/3, 2/3;1; 1 + t^{-3})-\overline{\rho}\cdot {_2F_1}(1/3, 2/3;1; -t^{-3})}{ {_2F_1}(1/3, 2/3;1; 1 + t^{-3}) +\rho\cdot {_2F_1}(1/3, 2/3;1; -t^{-3})}.
\end{equation}
Let $\tau$ be as in the statement of Proposition~\ref{PiAGM_rho}. Hence, if we let $\gamma={-t}^{-3}$, then we obtain
\begin{equation}\label{case1}
c_{\rho}(t)= \left(\frac{\Gamma(1 / 3)^3}{3\Gamma(2 / 3)^3}\right) \left(\frac{{\color{black}-}\sqrt{3}i \tau -{\color{black}\overline{\rho}} }{{\color{black}-}\sqrt{3}i \tau + {\color{black}\rho}}\right)
= \left(\frac{\Gamma(1 / 3)^3}{3\Gamma(2 / 3)^3}\right)\cdot \frac{(3\tau-1)-\rho}{(3\tau-1)-\overline{\rho}}.
\end{equation}
and {\color{black} by direct calculation we find that} $s_{\rho}^{-1}(c_{\rho}(t)){\color{black}=3\tau-1}. $ The proof of the theorem in this case follows from Proposition \ref{PiAGM_rho}, which gives
$$j(\tau) = 27t^3 \left(\frac{8 - t^3}{1 + t^3}\right)^3.$$

It turns out that the same method of proof works for all possible $0<|t|<1.$ One merely needs to keep track of  the branch cut crossings which arise in the hypergeometric transformation formulae.
There are six possible cases for $\arg(t)$, and the method applies {\it mutatis mutandis.}   
Apart from the Gamma-ratios, one obtains six different expressions for $c_{\rho}(t)=h_{\rho}(t)/g_{\rho}(t)$  in (\ref{crhoformula}),  where $h_{\rho}(t)$ and $g_{\rho}(t)$ are always of the form
$$
t^{-1}\left(\rho^a {}_2F_1(1/3,2/3;1;1+t^{-3})+\rho^b{}_2F_1(1/3,2/3;1;-t^{-3})\right).
$$
The exponents $a$ and $b$ depend on $\arg(t),$ and they are recorded in
Table 1. Moreover, the table gives the resulting  formula for 
$$
\widehat{c}_{\rho}(t):=c_\rho(t)/2\pi \Omega_{\rho}^2,
$$ 
which is the last factor of (\ref{case1}) when $\arg(t)\in [0,\pi/3).$
To complete the proof, we make use of the modular transformation properties of the $j$-function,
 and we note 
 that the corresponding
$(s_{\rho}^{-1}(c_{\rho}
(t))+1)/3$ values are indeed $\SL_2(\Z)$-equivalent.

\begin{center}
\begin{table}[H]
\begin{tabular}{|r|c|c|c|c|c|c|}
\hline
$\arg(t)$ & $[0,\pi/3)$ & $(\pi/3,2\pi/3)$ & $(2\pi/3,\pi)$ & $(\pi,4\pi/3)$ & $(4\pi/3,5\pi/3)$ & $(5\pi/3,2\pi)$ \\
\hline
\hline
$g_\rho(t)$ & $a=0,b=1$ & $a=2,b=1$ & $a=2,b=3$ & $a=4,b=3$ & $a=4,b=5$ & $a=0,b=5$ \\
\hline
$h_\rho(t)$ & $a=0,b=2$ & $a=4,b=2$ & $a=4,b=0$ & $a=2,b=0$ & $a=2,b=4$ & $a=0,b=4$ \\
\hline
$\widehat{c}_\rho(t)$ & $\frac{(3\tau-1)-\rho}{(3\tau-1)-\overline{\rho}}$ & $\frac{-1/(3\tau+1)-\rho}{-1/(3\tau+1)-\overline{\rho}}$ & $\frac{-1/(3\tau-2)-\rho}{-1/(3\tau-2)-\overline{\rho}}$ & $\frac{(3\tau+1)/(3\tau+2)-\rho}{(3\tau+1)/(3\tau+2)-\overline{\rho}}$ & $\frac{(3\tau-2)/(3\tau-1)-\rho}{(3\tau-2)/(3\tau-1)-\overline{\rho}}$ & $\frac{(3\tau+2)-\rho}{(3\tau+2)-\overline{\rho}}$\\
\hline
\end{tabular}
\vspace{5pt}
\caption{ }
\end{table}
\end{center}

\noindent
{\bf Case of $i$:}  The case where $t=0$ is confirmed by direct calculation. Therefore, we assume
  that $0<|t|<1/2,$ which ensures
 the convergence of the relevant functions.
We begin by noting that
\begin{displaymath}
\begin{split}
_2F_1\left(1/4, 1/4; 1/2; 4t^2\right)&=\sum_{n=0}^{\infty}\frac{(4t^2)^n}{n!(1/2)(3/2)\cdots ((2n-1)/2)}\cdot \left(\frac14\cdot \frac54 \cdots \frac{4n-3}{4}\right)^2\\
&= \sum_{n=0}^{\infty}\frac{t^{2n}}{(2n)!}\prod_{m=1}^n(4m-3)^2\\
&=g_i(t).
\end{split}
\end{displaymath}
Applying Proposition~\ref{hyper2} (1) with $a=\frac14, b=\frac14,$ and $z=4t^2,$  we obtain
\begin{equation}\label{i1}
g_i(t)=\frac{\Gamma(3/4)^2}{2\sqrt{\pi}}\left(\  _2F_1\left(1/2, 1/2; 1; \frac{1-2t}{2}\right) + {_2F_1\left (1/2, 1/2; 1; \frac{1+2t}{2}\right)} \right).
\end{equation}
Similarly, we have 
\begin{displaymath}
\begin{split}
t\cdot {_2F_1\left(3/4, 3/4; 3/2; 4t^2\right)}&= t\cdot \sum_{n=0}^{\infty}\frac{(4t^2)^n}{n!(3/2)(5/2)\cdots ((2n+1)/2)}\cdot \left( \frac34\cdot \frac74\cdots \frac{4n-1}{4} \right)^2\\
&=\sum_{n=0}^{\infty} \frac{t^{2n+1}}{(2n+1)!}\prod_{m=1}^n(4m-1)^2 \\
&=h_i(t).
\end{split}
\end{displaymath}
Applying Proposition~\ref{hyper2} (2) with $a=\frac34, b=\frac34,$ and $z=4t^2$ gives
\begin{equation}\label{i2}
h_i(t)= -\frac{\Gamma(1/4)^2}{8\sqrt{\pi}} \left( _2F_1\left(1/2, 1/2; 1; \frac{1-2t}{2}\right) -{_2F_1\left(1/2, 1/2; 1; \frac{1+2t}{2} \right)} \right) .
\end{equation}
Therefore, by combining (\ref{i1}) and (\ref{i2}) we Proposition~\ref{Flat} (2), we find that
$$
c_i(t)=\left(\frac{\Gamma(1/4)^2}{-4\Gamma(3/4)^2}\right) \left(\frac{_2F_1(1/2, 1/2; 1; \frac{1-2t}{2}) - {_2F_1(1/2, 1/2; 1; \frac{1+2t}{2})} }{ _2F_1(1/2, 1/2; 1; \frac{1-2t}{2}) + {_2F_1(1/2, 1/2; 1; \frac{1+2t}{2})}}\right).
$$
Finally, thanks to Proposition~\ref{PiAGM_i}, if we let $\lambda =\frac{1-2t}{2} $, then we obtain
$$
c_i(t)=\left(\frac{\Gamma(1/4)^2}{4\Gamma(3/4)^2}\right)\left( \frac{-i+\tau}{i+\tau}\right),
$$
{\color{black}and by direct calculation we find that $s_{i}^{-1}(c_{i}(t))=\tau. $ Hence, we have}
$$
j(\tau)=64\cdot \frac{(3+4t^2)^3}{(1-4t^2)^2},
$$
and this completes the proof.

\section{Examples}\label{examples}
Here we offer numerical examples of the main theorem.

\begin{example}
Choosing $t=\sqrt 3-1$ in the $\tau_*=\rho$ case gives the value
$$
27t^3\left(\frac{8-t^3}{1+t^3}\right)^3=1728.
$$
Using the first 1000 terms of $c_\rho(t)$, we find that $c_\rho(\sqrt{3}-1)=0.691592015\dots,$ and we obtain
$$
\frac{s_\rho^{-1}(c_\rho(t))+1}{3}= (1+i)\cdot 0.500000000\dots.
$$
We note that the classical theory of $_2F_1$-Gaussian hypergeometric functions, applied to
(\ref{rho1}) and (\ref{rho2}), gives the exact value $(1+i)/2.$
This corresponds to the famous singular modulus
$$j((1+i)/2)=1728.$$
\end{example}

\begin{example}
The $\tau_*=\rho$ case of Theorem~\ref{main} when $t=1/2$ involves the $j$-invariantis exactly
$$
27t^3\left(\frac{8-t^3}{1+t^3}\right)^3=\frac{9261}{8}=1157.625.
$$
Using the first 1000 terms of $c_\rho(t)$, we find that
$c_{\rho}(1/2)=0.490175\dots,$
which in turn gives
$$
\frac{s_\rho^{-1}(c_\rho(1/2))+1}{3} \approx 0.50000000+0.424026095i
$$
Numerically, we happily find that
$$
j\left(\frac{s_\rho^{-1}(c_\rho(1/2))+1}{3}\right)\approx j(0.50000000+0.424026095i)=1157.625000\dots.
$$
\end{example}

\begin{example}
The $\tau_*=i$ case of Theorem~\ref{main} when $t=1/3$ involves the $j$-invariant
$$
64\cdot  \frac{(3+4t^2)^3}{(1-4t^2)^2}=\frac{1906624}{225}=8473.884444\dots.
$$
Using the first 1000 terms of $c_i(t)$, we find that
$c_{i}(1/3)=0.3830612321\dots,$ which in turn gives
$$
s_i^{-1}(c_i(1/3))\approx s_{i}^{-1}(0.3830612321)\approx 1.4243556206i.
$$
Numerically, we happily find that
$$
j\left(s_i^{-1}(c_i(1/3))\right)\approx  j(1.4243556206i)\approx 8473.884444\dots.
$$
\end{example}

\end{document}